\documentclass[12pt]{amsart}
\usepackage{a4wide}
\usepackage{amsfonts,amssymb,amsmath,color}

\usepackage{amssymb}
\usepackage{pdfsync}
\usepackage{wrapfig,epsfig}
\usepackage{color}
\usepackage{hyperref}

\newtheorem{definition}{Definition}[section]
\newtheorem{corollary}[definition]{Corollary}

\newtheorem{proposition}[definition]{Proposition}
\newtheorem{theorem}[definition]{Theorem}

\makeatletter
\@addtoreset{equation}{section}

\newcommand{\bes}{\begin{displaymath}}
\newcommand{\ees}{\end{displaymath}}
\newcommand{\be}{\begin{equation}}
\newcommand{\ee}{\end{equation}}
\newcommand{\ba}{\begin{eqnarray}}
\newcommand{\ea}{\end{eqnarray}}
\newcommand{\bas}{\begin{eqnarray*}}
\newcommand{\eas}{\end{eqnarray*}}
\newcommand{\@Bbb}[1]{\ensuremath{\Bbb #1}}

\newcommand{\B}{{\@Bbb B}}
\newcommand{\C}{{\@Bbb C}}
\newcommand{\E}{{\@Bbb E}}
\newcommand{\F}{{\@Bbb F}}
\newcommand{\G}{{\@Bbb G}}
\renewcommand{\P}{{\@Bbb P}}
\newcommand{\bbP}{{\P}}
\newcommand{\Q}{{\@Bbb Q}}
\newcommand{\bQ}{{\@Bbb Q}}
\newcommand{\N}{{\@Bbb N}}
\newcommand{\R}{{\@Bbb R}}
\newcommand{\T}{{\@Bbb T}}
\newcommand{\bbR}{{\@Bbb R}}
\newcommand{\W}{{\@Bbb W}}
\newcommand{\Z}{{\@Bbb Z}}
\newcommand{\bbZ}{{\@Bbb Z}}

\newcommand{\si}{\sigma}

\newcommand{\om}{\omega}

\newcommand{\ep}{\varepsilon}

\newcommand{\mbR}{\mathbb R}

\newcommand{\@s}[1]{\ensuremath{\mathcal #1}}
\newcommand{\cA}{\@s A}
\newcommand{\cB}{\@s B}
\newcommand{\cC}{\@s C}
\newcommand{\cD}{\@s D}
\newcommand{\cE}{\@s E}
\newcommand{\cF}{\@s F}
\newcommand{\cG}{\@s G}
\newcommand{\cH}{\@s H}
\newcommand{\cI}{\@s I}
\newcommand{\cJ}{\@s J}
\newcommand{\cal}{\mathcal}
\newcommand{\cK}{\@s K}
\newcommand{\cL}{\@s L}
\newcommand{\cN}{\@s N}
\newcommand{\cM}{\@s M}
\newcommand{\cO}{\@s O}
\newcommand{\cP}{\@s P}
\newcommand{\cQ}{\@s Q}
\newcommand{\cR}{\@s R}
\newcommand{\cS}{\@s S}
\newcommand{\cT}{\@s T}
\newcommand{\cU}{\@s U}
\newcommand{\cV}{\@s V}
\newcommand{\cW}{\@s W}
\newcommand{\cX}{\@s X}
\newcommand{\cY}{\@s Y}
\newcommand{\cZ}{\@s Z}

\def\qed{\mbox{$\square$}}
\newcommand{\@bm}[1]{\ensuremath{\mathbf #1}}
\newcommand{\bma}{\@bm a}\newcommand{\bmA}{\@bm A}
\newcommand{\bmb}{\@bm b}\newcommand{\bmB}{\@bm B}
\newcommand{\bmc}{\@bm c}\newcommand{\bmC}{\@bm C}
\newcommand{\bmd}{\@bm d}\newcommand{\bmD}{\@bm D}
\newcommand{\bme}{\@bm e}
\newcommand{\bmf}{\@bm f}\newcommand{\bmF}{\@bm F}
\newcommand{\bmg}{\@bm g}\newcommand{\bmG}{\@bm G}
\newcommand{\bmh}{\@bm h}\newcommand{\bmH}{\@bm H}
\newcommand{\bmi}{\@bm i}\newcommand{\bmI}{\@bm I}
\newcommand{\bmj}{\@bm j}
\newcommand{\bmk}{\@bm k}\newcommand{\bmK}{\@bm K}
\newcommand{\bml}{\@bm l}
\newcommand{\bmm}{\@bm m}\newcommand{\bmM}{\@bm M}
\newcommand{\bmn}{\@bm n}
\newcommand{\bmo}{\@bm o}
\newcommand{\bmp}{\@bm p}
\newcommand{\bmq}{\@bm q}\newcommand{\bmQ}{\@bm Q}
\newcommand{\bmr}{\@bm r}
\newcommand{\bms}{\@bm s}\newcommand{\bmS}{\@bm S}
\newcommand{\bmt}{\@bm t}
\newcommand{\bmu}{\@bm u}\newcommand{\bmU}{\@bm U}
\newcommand{\bmw}{\@bm w}\newcommand{\bmW}{\@bm W}
\newcommand{\bmv}{\@bm v}\newcommand{\bmV}{\@bm V}
\newcommand{\bmx}{\@bm x}\newcommand{\bmX}{\@bm X}
\newcommand{\bx}{x}
\newcommand{\bmy}{\@bm y}\newcommand{\bmY}{\@bm Y}\newcommand{\by}{\@bm y}
\newcommand{\bmz}{\@bm z}\newcommand{\bmZ}{\@bm Z}

\newcommand{\bbE}{\mathbb E}

\newcommand{\bmzero}{\@bm 0}

\newcommand{\ga}{\gamma}

\newcommand{\@g}[1]{\ensuremath{\mathfrak #1}}
\newcommand{\gA}{\@g A}
\newcommand{\gD}{\@g D}
\newcommand{\gJ}{\@g J}
\newcommand{\gF}{\@g F}
\newcommand{\gM}{\@g M}
\newcommand{\gR}{\@g R}

\newcommand{\commentout}[1]{{}}

\makeatother
\begin{document}

\date{\today}

\title[Estimates for a divergence form operator]{A probabilistic proof of apriori $L^p$ estimates
 for   a class of  divergence form elliptic operators}

\author{Tymoteusz Chojecki}
\address{Tymoteusz Chojecki: Institute of Mathematics,  UMCS
  \\ pl. Marii Curie-Sk\l odowskiej 1 \\ 20-031,
Lublin, Poland.} 
\email{{\tt chojecki.tymoteusz@gmail.com}}

\author{Tomasz Komorowski}
\address{Tomasz Komorowski: Institute of Mathematics, Polish Academy
  Of Sciences\\ ul. \'Sniadeckich 8\\ 00-636 Warsaw, Poland.} 
\email{{\tt komorow@hektor.umcs.lublin.pl}}

\maketitle

\begin{abstract}
Suppose that ${\cal L}$ is  a divergence form differential operator of the
form
${\cal
  L}f:=(1/2) e^{U}\nabla_x\cdot\big[e^{-U}(I+H)\nabla_x
  f\big]$, where $U$ is scalar valued, $I$ identity matrix and $H$  an
anti-symmetric matrix valued function. The coefficients are not assumed to be bounded, but are
$C^2$ regular. 
We show that  if $Z=\int_{\bbR^d}e^{-U(x) }dx<+\infty$ and   the supremum of the
numerical range of matrix $-\frac12\nabla^2_x U+\frac12\nabla_x\left\{\nabla_x\cdot H-[\nabla_x U]^TH\right\}$
satisfies some exponential integrability condition with respect to
measure $d\mu=Z^{-1}e^{-U}dx$,  then
for any $1 \le p<q<+\infty$ there exists a constant $C>0$ such that  $\left\| f\right\|_{W^{2,p}(\mu)}\le
C\Big(\left\|{\cal
    L}f\right\|_{L^q(\mu)}+\left\|f\right\|_{L^q(\mu)}\Big)$ for
$f\in C_0^\infty(\bbR^d)$. Here
   $W^{2,p}(\mu)$ is the Sobolev space of functions that are $L^p(\mu)$ integrable with two  derivatives. Our proof is probabilistic and relies on an application of the Malliavin calculus.

\end{abstract}

\section{Introduction}

In the present note we formulate some apriori bound for a divergence form differential operator of the
form
\begin{align}\label{CL}
&{\cal
  L}f(x):=\frac{1}{2}e^{U(x)}\nabla_x\cdot\Big[e^{-U(x)}(I+H(x))\nabla_x
  f(x)\Big]\nonumber\\
&
=\sum_{i,j=1}^d \frac{1}{2}e^{U(x)}\partial_{x_i}\left\{e^{-U(x)}\left[\left(\delta_{i,j}+H_{i,j}(x)\right)\partial_{x_j}f(x)\right]\right\}, \quad f\in C_0^{\infty}(\bbR^d),\,x\in\bbR^d.
\end{align}
Here $U:\bbR^d\to \bbR$
and $H(x)=[H_{i,j}(x)]$ is a $d\times d$ antisymmetric matrix valued
function, i.e.
\begin{equation}
\label{010904-19}
H_{i,j}(x)=-H_{j,i}(x),\quad x\in\bbR^d,\,i,j=1,\ldots,d.
\end{equation}
We do not assume that the potential $ U(x)$ and antisymmetric part
$ H(x)$ are bounded. Instead we assume that
\begin{equation}\label{stalaZ}
Z:=\int_{\bbR^d}e^{-U(x) }dx<+\infty.
\end{equation} 
and furthermore they satisfy some exponential integrability condition with respect to the measure 
\begin{equation}\label{miara}
\mu(dx):=\dfrac{1}{Z}e^{-U(x) }dx,
\end{equation}
see \eqref{zalU} below.

In our main result, see Theorem \ref{TwProf} and Corollary \ref{Uwaga} below,  we show that for any $q>p>1$, there exists
constant $C>0$, depending only on $p,q$, dimension $d$ and exponential moments of the coefficients, with
respect to $\mu$,  such that
\begin{equation}\label{Proff11zz}
\|f\|_{W^{m,p}(\mu)}\le C(\|f\|_{L^q(\mu)}+\|{\cal
  L}f\|_{L^q(\mu)}), \quad f\in
C_0^\infty(\bbR^d). 
\end{equation} 
Here, for any integer $m\ge1$ and $p\in[1,+\infty)$, we denote by 
$W^{m,p}(\mu)$ the  respective Sobolev space, defined  as the completion of $C_0^\infty(\bbR^d)$ in the norm
\begin{equation}
\label{sobolev}
\|f\|_{W^{m,p}(\mu)}:=\left\{\|f\|_{L^{p}(\mu)}^p+\sum_{k=1}^m\|\nabla^k_xf\|_{L^{p}(\mu)}^p\right\}^{1/p},\quad f\in C_0^\infty(\bbR^d).
\end{equation}
These type of estimates have been used by the authors in proving  homogenization of advection equations with non-stationary, locally ergodic coefficients, see \cite{chk}.

A priori estimates  for Poisson equation in case of Sobolev norms with
respect to the flat Lebesgue measure, and the Laplacian, or a uniform elliptic operator with bounded coefficients, are  classical, see e.g. \cite{Gilbarg}. 
For ${\cal L}$ that is an Ornstein-Uhlenbeck operator (then $U$ and $H$ are quadratic polynomials, and  $\mu$ is therefore Gaussian) one can in fact assume that $p=q$ in \eqref{Proff11zz}, see \cite{MP1}. This result has been generalised in \cite{MP2} to cover some generalisation of Ornstein-Uhlenbeck operators, quite different from the one considered in the present paper. Some further results concerning  global $L^p$ bounds for divergence form operators    with lower oreder coefficients in an appropriate Morrey class can be found in \cite{Transirico}.



 Our proof of \eqref{Proff11zz} is probabilistic and relies on an application of the Malliavin calculus.
 The gradient of function $f$ can be represented using the Fr\'{e}chet derivative of the
stochastic flow corresponding to the diffusion with the generator ${\cal L}$,  see \eqref{mamy} and \eqref{010706d} below. The principal
novelty of the paper is to replace the Fr\'{e}chet by the Malliavin
derivative  in the direction of a random path that is 
adapted with respect to the natural  filtration of the flow and such
that, after some initial boundary layer, the respective   Fr\'{e}chet and Malliavin
derivatives coincide, see \eqref{wzor1} and \eqref{010604-19x}
below. Finding such a path requires to solve a simple linear control
problem, see Section \ref{sec6.4}. The solution also yields an
estimate of the moments of the random path, see \eqref{030707}. {We mention here that an analogous  argument has been used by D. Bell in \cite{bell}
in the context of finding a representation of a certain class of vector fields, defined over the space of paths on a Riemannian manifold, in terms of an appropriately defined divergence operator of the field.}

The organisation of the paper is as follows: the basic definitions and
assumptions are stated in Section \ref{basic}, where also the main
results of the paper are rigorously stated, see Section \ref{MR}. The proofs of
the results are given throughout Section \ref{sec3}. Finally, in the
appendix we give the proofs of the auxiliary results, mainly from
the theory of diffusions that are used throughout the  paper.

\bigskip

{

\subsection*{Acknowledgements}
Both authors acknowledge the support of the National Science Centre:
NCN grant 2016/23/B/ST1/00492. T.K. wishes to express his gratitude to Professors D. Bakry and
K. Oleszkiewicz for useful discussions on the topic of the paper.


\section{Preliminaries and formulation of  the main result}
\bigskip

\label{basic}

We 
assume that the coefficients of ${\cal L}$ are sufficiently smooth, so there is no
issue with the definition of the stochastic flow
$\left(X(t,x)\right)_{t\ge0,x\in\bbR^d}$ and semigroup $(P_t)_{t\ge0}$
corresponding to the operator, see Section \ref{sec3.1} below.

\subsection{Differential operator} 

\label{sec-diff}

\begin{itemize}
\item[{\bf A1)}] We assume that  $U\in C^2(\mbR^d)$ and 
$H=[H_{i,j}]$ is a $d\times d$ matrix valued
function with $C^2$ smooth entries that  satisfy \eqref{010904-19}.

\end{itemize}

We can rewrite $\mathcal{L}$ in the following way
\begin{equation}\label{calL}
\mathcal{L}f=Lf+\mathcal{A}f,\quad f\in C_0^{\infty}(\bbR^d),
\end{equation}
where operator $L$ is defined as
\begin{equation}
\label{L}
Lf:=\frac{1}{2}e^{U}\nabla_x\cdot\left(e^{-U}\nabla_x
  f\right)
=\frac{1}{2}\Delta_x f-\frac12\nabla_x U\cdot \nabla_x f,\quad f\in C_0^{\infty}(\bbR^d)
\end{equation}
and
\begin{equation}\label{CA}
{\cal A}f:=\frac12b\cdot\nabla_x f,\quad f\in C_0^{\infty}(\bbR^d).
\end{equation}
where $b^T=[b_1,\ldots,b_d]$, is given by
\begin{equation}\label{CA1}
b:=\nabla_x\cdot H-[\nabla_x U]^TH,
\end{equation}
or coordinatewise
\begin{equation}\label{warB2}
b_j(x):=\sum_{i=1}^d\left(\partial_{x_i}H_{i,j}(x)-\partial_{x_i}U(x)H_{i,j}(x)\right),\quad x\in\mbR^d,\,j=1,\ldots,d.
\end{equation}

\subsection{Main results}

\label{MR}

To ensure global existence of the stochastic flow
corresponding to ${\cal L}$, see Section \ref{sec3.1} below, we assume that:
\begin{itemize}
\item[{\bf A2)}] there exists $\rho >0$ such that
\begin{equation}
\label{011104-19}
\int_{\bbR^d}|x|^{2\rho}\mu(dx)+\int_{\bbR^d}(|x|^2+1)^{\rho-1}\left|
  x\cdot\left[-\nabla_x U(x)+b(x)\right]\right| \mu(dx)<+\infty.
\end{equation}
\end{itemize}
Let
    ${\cal U} $ be the matrix  valued function, given by
\begin{equation}
\label{cU}
{\cal U}:=-\frac12\nabla^2_x U+\frac12\nabla_x b=-\frac12\nabla^2_x U+\frac12\nabla_x\left\{\nabla_x\cdot H-[\nabla_x U]^T
  H\right\},
\end{equation}
or entrywise ${\cal
  U}=[{\cal
  U}_{j,j'}]$, with 
$$
{\cal
  U}_{j,j'}:=-\frac12\partial^2_{x_j,x_{j'}}U+\frac12\sum_{i=1}^d\partial_{x_{j'}}\left\{\partial_{x_i}H_{i,j}-H_{i,j}\partial_{x_i}U\right\},\quad
j,j'=1,\ldots,d.
$$
We let 
$$
{\frak u}:=\max_{|\ell|=1}\langle {\cal U}\ell,\ell\rangle_{\bbR^d}.
$$
\begin{itemize}
\item[\textbf{A3)}]{\bf Integrability condition.} Assume that there exists $\ga_0>0$ such that
\begin{equation}\label{zalU}
E(\ga_0):=\int_{\mbR^d} \exp\left\{\ga_0 {\frak
    u}\vee 0\right\}d\mu<+\infty.
\end{equation}
\end{itemize}
Here $a\vee b:=\max\{a,b\}$.

Our  main result gives an estimate of the $L^p(\mu)$ norm of the
  gradient of $f\in C_0^\infty(\bbR^d)$ in terms of the $L^q(\mu)$
  norms of $f$ and ${\cal L}f$, with $q>p$.
{
\begin{theorem}\label{TwProf}
Suppose that \eqref{stalaZ} and assumptions $A 1)-A 3)$ hold. Then, for any $p\in[1,+\infty)$
and  $q>p$ there exists $C(d,r,\ga_0)>0$, depending only on
the dimension $d$, $\ga_0$ and $r$, given by $1/p=1/q+1/r$ , for which
\begin{equation}\label{Proff1aa}
\left\|\nabla_x f\right\|_{L^{p}(\mu)}\le
C(d,r,E(\ga_0)) \left(\left\|{\cal L}f\right\|_{L^q(\mu)}+\left\|f\right\|_{L^q(\mu)}\right),\quad f\in C_0^\infty(\bbR^d).
\end{equation}
\end{theorem}}
 The proof of this theorem is contained in Section \ref{Dowod1}.

Using an
 argument based on Calderon-Zygmund estimates we can also conclude the
 following.
\begin{corollary}\label{Uwaga}
 Suppose that all the hypotheses of  Theorem $\ref{TwProf}$ remain in force. 
In addition, we assume  that
\begin{equation}
\label{zalozenie}
\ell_{2,*}:=\sum_{j,j'=1}^d\|H_{j,j'}\|_{W^{1,r}(\mu)}+\|U\|_{W^{2,r}(\mu)}<+\infty,\quad \forall r\in[1,+\infty).
\end{equation}
Then, for any $q>1$
and $p\in[1,q)$ there exists $C(d,p,q,E(\ga_0),\ell_{2,*})>0$, depending only on the
indicated parameters, such
that
\begin{equation}\label{Proff1}
\left\|\nabla_x^2 f\right\|_{L^{p}(\mu)}\le
C\left(\left\|{\cal L}f\right\|_{L^q(\mu)}+\left\|f\right\|_{L^q(\mu)}\right),\quad f\in C_0^\infty(\bbR^d).
\end{equation}
\end{corollary}
The proof of the corollary is presented in Section \ref{dowodUw}.

\section{Proofs of the Main Results}
\label{sec3}

\subsection{Diffusion corresponding to ${\cal L}$} 

\label{sec3.1}

It follows from the assumptions made about the coeficients of ${\cal L}$
that they are locally Lipschitz,
therefore we can define  a
diffusion $\left(X(t;x)\right)_{0\le t<\tilde{\frak e}_x}$ via
an It\^o stochastic differential equation 
\begin{equation}\label{xt}
\begin{aligned}
&dX(t;x)=\left\{-\frac12\nabla_x U(X(t;x))+\frac12b(X(t;x))\right\}dt
+ dw(t), \quad 0\le t<\tilde{\frak e}_x
&X(0;x)=x,
\end{aligned}
\end{equation}
 where 
 $\left(w(t)=(w_1(t),\ldots,w_d(t))\right)_{t\ge0}$ is standard $d-$dimensional Brownian
 motion and $\tilde{\frak e}_x$ is an explosion time, see
 \eqref{frakex} below. The solutions satisfy the flow property, i.e.
$$
X(t;X(s,x))=X(t+s;x),\quad t+s<\tilde{\frak e}_x,\,t,s\ge0,\,x\in\bbR^d
$$
for a.s. realization of the field $\left(X(t;x)\right)_{0\le
  t<\tilde{\frak e}_x,\,x\in\bbR^d}$ , see e.g. Theorem 4.2.5, p. 126
of \cite{kunita}.


Integrating by parts we conclude  that $L$ is symmetric and ${\cal A}$
is anti-symmetric  with respect to $\mu$,
 i.e.
\begin{equation}
\label{sym}
\int_{\bbR^d}Lfgd\mu=\int_{\bbR^d}fLgd\mu
\end{equation}
and
\begin{equation}
\label{e2}
\int_{\mbR^d}{\cal A}fgd\mu=-\int_{\mbR^d}f{\cal A}gd\mu,\quad f,g\in C_0^{\infty}(\mbR^d).
\end{equation}

The processes $\left(X(t;x)\right)_{t\ge0}$ could, in principle,
explode in finite time. This scenario is however precluded, thanks to
condition A2). More precisely,
for any path $\pi:[0,+\infty)\to \bbR^d$ and $R>0$ define
\begin{equation}
\label{tauR}
\tau_R(\pi):=\inf_{t\ge0}[|\pi(t)|\ge R],
\end{equation}
with the convention that $\tau_R(\pi)=+\infty$, if the set over which
we take infimum is empty.

The function $R\mapsto\tau_R $ is increasing. We can define therefore
$$
{\frak e}(\pi):=\lim_{R\to+\infty}\tau_R(\pi).
$$
Suppose that the random vector $\eta$ is distributed according to
$\mu$ and  independent of  $\left(w(t)\right)_{t\ge0}$.  Consider the process
$$
X(t):=X(t;\eta),\quad 0\le t<\tilde{\frak e}:={\frak e}\left(X(\cdot;\eta)\right).
$$
\begin{proposition}
\label{prop010804-19}
Suppose that \eqref{stalaZ} and conditions A1) - A2) hold.
Then,
\begin{equation}
\label{fe}
\bbP[\tilde{\frak e}<+\infty]=0.
\end{equation}
\end{proposition}

As a result of the above proposition, we can define the process
$\left(X(t)\right)_{t\ge0}$ for all times.
\begin{proposition}
\label{prop020804-19}
The process $\left(X(t)\right)_{t\ge0}$ is stationary, i.e. for any
$0\le t_1<t_2<\ldots<t_n$, bounded measurable functions $f_1,\ldots,f_n$ and $h\ge0$ we have
\begin{equation}
\label{stat}
\bbE\left[\prod_{i=1}^nf_i(X(t_i+h))\right]=\bbE\left[\prod_{i=1}^nf_i(X(t_i))\right].
\end{equation}
\end{proposition}
The proofs of  Propositions \ref{prop010804-19} and
\ref{prop020804-19} are contained in Section \ref{sec2.1} of Appendix.

It turns out that we can define the diffusion starting from any $x$
for all times. Let 
\begin{equation}
\label{frakex}
\tilde{\frak e}_x:={\frak e}\left(X(\cdot;x)\right).
\end{equation}
\begin{proposition}
\label{prop030804-19}
For any $x\in\bbR^d$ we have
\begin{equation}
\label{fet}
\tilde{\frak e}_x<+\infty,\quad\mbox{a.s.}
\end{equation}
\end{proposition}
The proof of this result is contained in Section \ref{sec2.3} of Appendix. Thanks
to the above result the trajectories $t\mapsto X(t;x)$ are defined for
all $t\ge0$, a.s. for any $x\in\bbR^d$.
 


Let $P_t(x,dy)$ be the transition probabilities corresponding to
$\left(X(t;x)\right)_{t\ge0}$. Thanks to Proposition \ref{prop020804-19} we
have
\begin{equation}
\label{statm}
\int_{\bbR^d}\mu(dx)P_t(x,A)=\mu(A),\quad A\in{\cal B}(\bbR^d),\,t>0.
\end{equation}
Define by  $(P_t)_{t\ge0}$ the respective transition probability semigroup, given by $P_0=I$,
\begin{equation}\label{polgrupa}
P_tf(x)=\bbE f(X(t;x))=\int_{\bbR^d}f(y)P_t(x,dy),\quad f\in B_b(\mbR^d), t>0,\,x\in\bbR^d.
\end{equation}
In consequence of  the path continuity of diffusions we conclude the following.
\begin{corollary}
\label{cor020804-19}
The semigroup $(P_t)_{t\ge0}$  extends to a $C_0$ semigroup on
$L^p(\mu)$ for any $p\ge1$.
\end{corollary}
We use an abbreviated notation $P:=P_1$.
 \begin{proposition}
\label{prop040804-19}
There exists a Borel measurable, strictly positive function
$p:\bbR^{2d}\to(0,+\infty)$ such that for every $f\in L^1(\mu)$
\begin{equation}
\label{PP}
Pf(x)=\int_{\bbR^d} p(x,y)f(y)\mu(dy),\quad \mu\,\mbox{a.e.}
\end{equation}
\end{proposition}
The proof of the proposition is contained in Section \ref{sec2.4} of Appendix.

As a consequence of \eqref{PP} and Theorem 6 of \cite{rudnicki} we get
\begin{corollary}
\label{cor011304-19}
For any $f\in L^2(\mu)$ such that
\begin{equation}
\label{050804-19}
\int_{\bbR^d}fd\mu=0
\end{equation}
 we have
\begin{equation}\label{Pt}
\lim_{t\to+\infty}\|P_tf\|_{L^2(\mu)}=0.
\end{equation}
\end{corollary}

\subsection{Proof of Theorem \ref{TwProf}}
\label{Dowod1}
 We assume with no loss of generality that $f\in C_0^\infty(\bbR^d)$
 satisfies \eqref{050804-19}.  We have
\begin{equation}\label{PG}
f(x)=-\int_0^{+\infty}P_t{\cal L}f(x)dt=-\int_0^{+\infty}\bbE[{\cal
  L}f(X(t;x))]dt,\quad x\in\bbR^d.
\end{equation}
Therefore, 
\begin{equation}\label{mamy}
f_{x_j}(x)=\int_0^{+\infty}v_j(t)dt,
\end{equation}
where $v_j(t)$ is given by 
\begin{align}
\label{010706d}
&
v_j(t):=-\partial_{x_j}P_t{\cal L}f(x)=-\bbE\left[(\nabla_x {\cal L}f)(X(t;x))\cdot \partial_{x_j}
X(t;x)\right]
\\
&
=-\sum_{i=1}^d\bbE\left[ ({\cal L}f)_{x_i}(X(t;x))
\xi_{i,j}(t)\right]\nonumber
\end{align}
and $\partial_{x_j}
X_i(t;x)=\xi_{i,j}(t)$ is  the Fr\'{e}chet derivative
of the stochastic flow $\bx\mapsto X(t;x)$.

 Differentiating \eqref{xt} with respect to the initial condition we conclude that
 \begin{align}\label{xi}\begin{aligned}
& \frac{d\xi_{i,j}}{dt}=\frac12\sum_{i'=1}^d\left[-\partial_{x_i,x_{i'}}^2U(X(t;x))+\partial_{x_{i'}}b_i(X(t;x))\right]\xi_{i',j},\\
&\xi_{i,j}(0)=\delta_{i,j},\ i,j=1,\ldots,d,
 \end{aligned}\end{align}
 where $\delta_{i,j}$ is the Kronecker symbol, i.e.  $\delta_{i,i}=1$ and $\delta_{i,j}=0$, if $i\not=j$.

For a given  random matrix valued field
$${\frak
  g}(t,{x}):=\left({\frak g}_j(t,{x})\right)_{j=1,\ldots,d}=\left[
  g_{i,j}(t,{x})\right]_{i,j=1,\ldots,d},\quad t\ge0,\ {x\in\mbR^d},
$$ 
and $r\ge 1$
we let
{\begin{align}
\label{030707x}
&\|\!|{\frak g}\|\!|_{r,t}:=\sum_{i,j=1}^d\left\{\int_{\bbR^{d}}\mu(dx)\bbE\left[\int_0^t|g_{i,j}{(s,x)}|^2ds\right]^{r/2}\right\}^{1/r}.
\end{align}
Let also
${\frak h}_j(t,{x})=(h_{1,j}(t,x),\ldots, h_{d,j}(t,{x}))$, $j=1,\ldots,d$, where
\begin{equation}\label{ggg}
h_{i,j}(t,{x}):=\int_0^tg_{i,j}(s,x)ds,\quad t\ge0,\ {x\in\mbR^d}.
\end{equation}
Suppose that $\left({\frak g}(t,\cdot)\right)_{t\ge0}$ is with respect to the natural filtration of
  $\left(w(t)\right)_{t\ge0}$. Treating the solution $X(t;x,w)$ of \eqref{xt} as the functional of the Wiener process $w(t)=(w_{i}(t))_{i=1,\ldots,d}$,
we define the Malliavin derivative of $X(t;x,w)$ in the direction ${\frak h}_j$
\begin{align*}
&\mathcal{D}_{\frak h_j}X(t;x,w):=\lim_{\ep\to0}\frac{1}{\ep}\left\{X(t;x,w+\ep
{\frak h_j})-X(t;x,w)\right\},
\end{align*}
 where the limit above is understood in the $L^2$ sense, see Def. 2.1, p. 35 of \cite{malliavin}.
 
  Denote
$$
\zeta_{i,j}(t):=\mathcal{D}_{{\frak h}_j}X_i(t;x,w),\quad i,j=1,\ldots,d
$$
 the components of the Malliavin derivatives in the directions ${\frak h}_1,\ldots,{\frak h}_d$. We can see, from \eqref{xt},
that they satisfy
\begin{align}\label{zetatheta}\begin{aligned}
& \frac{d\zeta_{i,j}}{dt}=\frac12\sum_{i'=1}^d\left[-\partial_{x_i,x_{i'}}U(X(t;x))+\partial_{x_{i'}}b_i(X(t;x))\right]\zeta_{i',j}+g_{i,j},\\
&\zeta_{i,j}(0)=0,\ i,j=1,\ldots,d,
\end{aligned}\end{align}
From the chain rule for the Malliavin
derivative, see Proposition 1.2.3, p. 28 of \cite{Nualart}, we obtain
\begin{align}
\label{020707vv}
\bbE\left[\mathcal{D}_{{\frak h}_j}{\cal L}f(X(t,x))\right]=\sum_{i=1}^d\bbE\left[\partial_{x_i}{\cal L}f( X(t,x))\zeta_{i,j}(t)\right].
\end{align}

The difference of the Fr\'echet and Malliavin derivatives
\begin{equation}\label{TH}
\Theta(t,x) := \partial_{x}X(t;x)-\mathcal{D}_{\frak h}X(t;x)=
(\Theta_{i,j}(t,x))_{i=1,\ldots,d}
\end{equation}
satisfies the following system of equations
\begin{align}\label{xitheta}\begin{aligned}
&
\frac{d\Theta_{i,j}(t,x)}{dt}=\frac12\sum_{i'=1}^d\left[-\partial_{x_i,x_{i'}}U(X(t;x))+\partial_{x_{i'}}b_j(X(t;x))\right]\Theta_{i',j}(t,x)-g_{i,j}(t,{x}),\\
&\Theta_{i,j}(0,x)=\delta_{i,j},\ i,j=1,\ldots,d.
\end{aligned}\end{align}
Accordingly, using \eqref{020707vv} together with \eqref{TH}, we get
\begin{align}
\label{020707}v_j(t)=\tilde v_j(t)-\bbE\left[\mathcal{D}_{{\frak h}_j}{\cal L}f( X(t;x))\right],
\end{align}
with $\tilde v_j(t)$ given by 
\begin{align}
\label{010706dv}
&\tilde 
v_j(t):
=-\sum_{i=1}^d\bbE\left[ ({\cal L}f)_{x_i}(X(t;x))
\Theta_{i,j}(t,x)\right].
\end{align}
Integrating by parts the second term on the right hand side of
\eqref{020707}, {see Theorem 2.1 p. 236 of \cite{bell} (and also Lemma 1.2.1 p. 25,
of \cite{Nualart})}, we
conclude that
\begin{align}\label{wzor1}\begin{aligned}
&v_j(t)=\tilde v_j(t)
-\sum_{i=1}^d\bbE\left[{\cal L}f(X(t;x))\int_0^tg_{i,j}(s,{x})dw_{i}(s)\right].
\end{aligned}\end{align}

We shall look for ${\frak g}^T(t,t_0,\bx)=\left[{\frak g}_1(t, t_0,{\bx}),\ldots,{\frak g}_d(t, t_0,{\bx})\right]$ - a column vector - that satisfies the following conditions:
\begin{itemize}
\item[i)] it is adapted with respect to the natural filtration of
  $\left(w(t)\right)_{t\ge0}$,
\item[ii)] given $t_0\in (0,t_*]$, where $t_*:=\gamma_0/r$ and parameter $\gamma_0$
as in $\eqref{zalU}$, we have  both  $\Theta(t,x)\equiv0$ and ${\frak
    g}(t, t_0,\bx)\equiv0$ for $t\ge t_0$,
  ${x} \in\bbR^{d}$,
\item[iii)] we have (cf \eqref{030707x})
{\begin{equation}
\label{030707}
\|\!|{\frak g}(\cdot;t_0)\|\!|_{r,t_0}<+\infty.
\end{equation}}
\end{itemize}

Suppose that we can construct such an object, that in what follows we call a control (we show how to do it
in Section \ref{sec6.4}). In the remaining part of the present section, we demonstrate how to conclude with its help the proof of \eqref{Proff1aa}.


Assume first that $t\ge
t_0$.  Thanks to ii), we conclude that then $\tilde v_j(t)\equiv0$.
Using formula \eqref{wzor1} and the Markov property of
$\left(X(t;x)\right)_{t\ge0}$ we can write
\begin{align}
\label{010604-19x}
&v_j(t)=-\sum_{i=1}^d\bbE\left[P_{t-t_0}{\cal L}f(X(t_0;x))\int_0^{t_0}g_{i,j}(s;t_0,x)dw_{i}(s)\right],\quad
\mbox{  for  }t\ge t_0,
\end{align}
which, invoking Corollary \ref{cor011304-19}, yields, upon integration from $t_0$ to $+\infty$,
\begin{align}
\label{010604-19a}
(P_{t_0}f)_{x_j}(x)=\int_{t_0}^{+\infty}v_j(t)dt=\sum_{i=1}^d\bbE\left[f( X(t_0;x))\int_0^{t_0}g_{i,j}(s;t_0,x)dw_{i}(s)\right].
\end{align}
First equality above comes from \eqref{PG}.
Denote ${\frak
    g}(t;t_0):={\frak
    g}(t; t_0,\eta)$, where $\eta$ is the stationary initial data vector for \eqref{xt}.

For $1/p=1/q+1/r$, we can write, by the
Burkholder-Davies-Gundy inequality
\begin{align}
\label{010604-19b}
&
\|\nabla_x P_{t_0}f\|_{L^p(\mu)}\le
 \|f\|_{L^q(\mu)}\left\{\bbE\left|\int_0^{t_0}{\frak
  g}(s;t_0)\cdot dw(s)\right|^r\right\}^{1/r}
\end{align}
\begin{align*}
&
\le  C(d,r)\|f\|_{L^q(\mu)}\left\{\bbE\left[\int_0^{t_0}{\rm tr}\,\left({\frak
  g}^T(s;t_0){\frak
  g}(s;t_0)\right)ds\right]^{r/2}\right\}^{1/r}
\\
&
= C(d,r)t_0^{1/2}\|f\|_{L^q(\mu)}\left\{\bbE\left[\frac{1}{t_0}\int_0^{t_0}{\rm tr}\,\left({\frak
  g}^T(s;t_0){\frak
  g}(s;t_0)\right)ds\right]^{r/2}\right\}^{1/r}\\
&
\le C(d,r)t_0^{1/2}\|f\|_{L^q(\mu)}\left\{\frac{1}{t_0}\int_0^{t_0}\bbE\left[{\rm tr}\,\left({\frak
  g}^T(s;t_0){\frak
  g}(s;t_0)\right)\right]^{r/2}ds\right\}^{1/r}.
\end{align*}
By virtue of estimate \eqref{060804-19} below we conclude that
\begin{align}
\label{010604-19bb}
&
\|\nabla_x P_{t_0}f\|_{L^p(\mu)}
\le
  C(d,r)t_0^{-1/2}\|f\|_{L^q(\mu)}\left\{\int_{\bbR^d}{\frak E}(t_0r{\frak u}(x))\mu(dx)\right\}^{1/r}, 
\end{align}
with ${\frak E}(\cdot)$ given by \eqref{E}. Therefore (see \eqref{010604-19a}), we obtain the following bound
\begin{align}\label{010708-18}
&\left\|\int_{t_0}^{+\infty}v_j(t)dt\right\|_{L^p(\mu)}\leq
  \frac{C}{t_0^{1/2}}\|f\|_{L^q(\mu)} ,\quad t_0\in(0,t_*)
\end{align}
and
\begin{equation}
\label{010410-19}
C:= C(d,r) \left\{\int_{\bbR^d}\exp\left\{\ga_0{\frak u}(x)\vee 0\right\}\mu(dx)\right\}^{1/r}
\end{equation}
It remains to consider the case when $t\in(0,t_0)$. We will use
estimate \eqref{010604-19bb}. Recall that $t_0$ appearing there can be
chosen arbitrarily from $(0,t_*]$. 
Applying \eqref{010604-19bb} with ${\cal L f}$ instead of $f$ and
$t\in(0,t_0]$ instead of $t_0$ we conclude that 
$$
\|v_j(t)\|_{L^p(\mu)}\le  \frac{C}{t^{1/2}}\|{\cal L}f\|_{L^q(\mu)}, \qquad t\in(0,t_0],
$$
with $C$ given by \eqref{010410-19}.
Therefore,
\begin{align}\label{010708-18aa}
&\left\|\int_0^{t_0}v_j(t)dt\right\|_{L^p(\mu)}\leq
  \int_0^{t_0}\|v_j(t)\|_{L^p(\mu)} dt\le C\|{\cal L}f\|_{L^q(\mu)}\int_0^{t_0}\frac{dt}{\sqrt{t}}=Ct_0^{1/2}\|{\cal L}f\|_{L^q(\mu)}.
\end{align}

From estimates \eqref{010708-18} and \eqref{010708-18aa} we conclude
that there exists a constant $C>0$ such that
\begin{equation}\label{Proff1b}
\left\|\nabla_x f\right\|_{L^{p}(\mu)}\le
C\left(t_0^{1/2}\left\|{\cal
      L}f\right\|_{L^q(\mu)}+t_0^{-1/2}\left\|f\right\|_{L^q(\mu)}\right),\quad
f\in C_0^\infty(\bbR^d),\quad t_0\in(0,t_*)
\end{equation}
and  \eqref{Proff1aa} follows.\qed

\subsection{Construction of a control ${\frak g}(t,t_0,x)$}

\label{sec6.4}

We construct first ${\frak g}(t,1,x)$. For a time being we suppress
writing arguments $1$ and $x$.
Denote by $C(t,s)=[C_{i,i'}(t,s)]_{i,i'=1,\ldots,d}$ the fundamental matrix of the system
\eqref{xi}.
It is a $d\times d$-matrix, which is the solution of the equation
\begin{equation}\label{eq16}
\frac{d}{dt}C(t,s)=A(t)C(t,s),\quad C(s,s)=I_{d},\quad t,s\ge0,
\end{equation}
where $I_{d}$ is the  identity $d\times d$-matrix and $A(t)=[A_{i,i'}(t)]_{i,i'=1,\ldots,d}$, where
\begin{equation}
\label{090707}
A_{i,i'}(t):=\frac12\left[-\partial_{x_i,x_{i'}}^2U(X(t,x))+\partial_{x_i}b_{i'}(X(t,x))\right].
\end{equation}

We  have 
$$
C(u,t)C(t,s)=C(u,s),\quad u,t,s\in\mbR.
$$
System \eqref{xitheta} can be rewritten as follows
\begin{equation}\label{eq15}
\frac{d\Theta(t)}{dt}=A(t)\Theta(t)-{\frak g}(t),\quad \Theta(0)=I_d,
\end{equation}
where $\Theta(t)=[\Theta_{ij}(t)]$ was defined in \eqref{TH} 
and $I_d$ is a $d$-dimensional identity matrix.

The solution of equation \eqref{eq15} can be expressed by the
Duhamel formula
\begin{equation}
\label{080707}
\Theta(t)=-\int_0^tC(t,s){\frak g}(s)ds+C(t,0),\quad t\ge0.
\end{equation}
We wish to show that $\Theta(t)\equiv 0$ for $t\ge1$, see condition ii).
For that purpose let
\begin{equation}
\label{100707}
{\frak g}(t):=C(t,0), \quad t\in[0,1],
\end{equation}
and let ${\frak g}(t):= 0$ for $t\ge 1$. The process  is adapted with
respect to the  natural filtration of $\left(w_t\right)_{t\ge0}$,
satisfying therefore condition 
i). Additionally, we have
$$
\Theta(1)=-\int_0^1C(1,s){\frak g}(s)ds+C(1,0)=-\int_0^1 C(1,0) ds+C(1,0)=0.
$$
Thus,  $\Theta(t)\equiv0$ for $t\ge1$. Condition ii) is therefore fulfilled.

In the general case we let
\begin{equation}
\label{gt}
{\frak g}(t;t_0,x):=\frac{1}{t_0}C(t,0),\quad t\in (0,t_0)\quad\mbox{and}\quad
{\frak g}(t;t_0,x):=0,\quad t\ge t_0.
\end{equation}
It remains to be checked that $\left({\frak g}(t,t_0,x)\right)_{t\ge0}$,
constructed above,
satisfies estimate \eqref{030707}.
From \eqref{eq16} we conclude that
\begin{equation}
\label{eq16a}
\frac{d}{dt}{\frak g}_i(t,t_0,x)=A(t) {\frak g}_i(t,t_0,x),\quad t\in[0,t_0],
\end{equation}
where $ {\frak g}_i(t,t_0,x)$ is the $i$-th column of  ${\frak g}(t,t_0,x)$.
Multiplying scalarly both sides of \eqref{eq16a} by ${\frak g}_i(t;t_0,x)$ we get
 $$
\frac{1}{2}\frac{d}{dt}|{\frak g}_i(t,t_0,x)|^2=A(t) {\frak
g}_i(t,t_0,x)\cdot {\frak g}_i(t;t_0,x)\le {\frak u}(X(t,x)) |{\frak g}_i(t,t_0,x)|^2.
$$
Thus,
by the Gronwall inequality we obtain
\begin{equation}
 \label{080708}
 |{\frak
   g}_i(t,t_0,x)|^2\leq\frac{1}{t_0^2}\exp\left\{2\int_0^{t} {\frak
     u}(X(s,x))ds\right\},\quad t\in[0,t_0]
 \end{equation}
Therefore, from \eqref{080708} (recall ${\frak
  g}(s;t_0)={\frak
  g}(s;t_0,\eta)$)
\begin{align*}
&\left\{\frac{1}{t_0}\int_0^{t_0}\bbE\left[{\rm tr}\,\left({\frak
  g}^T(s;t_0){\frak
  g}(s;t_0)\right)\right]^{r/2}ds\right\}^{1/r}\\
&
\le \sqrt{d} \left\{\frac{1}{t_0^{1+r}}\int_0^{t_0}\bbE\left[\exp\left\{r\int_0^{s} {\frak
     u}(X(\tau))d\tau\right\}\right]ds\right\}^{1/r}.
\end{align*}
Using the  Jensen inequality and then, subsequently, the stationarity
of $\left(X(t)\right)_{t\ge0}$,
 the right hand side estimates by
\begin{align*}
&
\sqrt{d}\left\{\frac{1}{t_0^{1+r}}\int_0^{t_0}\frac{ds}{s}\int_0^{s}\bbE\left[\exp\left\{r s{\frak
     u}(X(\tau))\right\}\right]d\tau\right\}^{1/r}
\\
&
= \sqrt{d} \left\{\frac{1}{t_0^{1+r}}\int_0^{t_0}\frac{ds}{s}\int_0^{s}d\tau\int_{\bbR^d}\exp\left\{r s{\frak
     u}(x)\right\}\mu(dx)\right\}^{1/r}\\
&
= \sqrt{d} \left\{\frac{1}{t_0^{1+r}}\int_0^{t_0}ds\int_{\bbR^d}\exp\left\{r s{\frak
     u}(x)\right\}\mu(dx)\right\}^{1/r}.
\end{align*}
Performing integration over the $s$ variable, we get that the utmost the right hand
side equals
\begin{align*}
&
\sqrt{ d} \left\{\frac{1}{t_0^{r}}\int_{\bbR^d}{\frak E}\left(r t_0{\frak
     u}(x)\right)\mu(dx)\right\}^{1/r},
\end{align*}
where ${\frak E}(0):=1$ and 
\begin{equation}
\label{E}
{\frak E}(x):=\frac{e^x-1}{x},\quad x\not=0.
\end{equation}
Note that 
$$
0<{\frak E}(x)\le
\left\{
\begin{array}{ll} e^x&\mbox{ for $x\ge 0$}\\
&
\\
 1\wedge |x|^{-1}&\mbox{ for $x<0$}.
\end{array}
\right.
$$
We have shown therefore that
\begin{align}
\label{060804-19}
&
 \left\{\frac{1}{t_0}\int_0^{t_0}\bbE\left[{\rm tr}\,\left({\frak
  g}^T(s;t_0){\frak
  g}(s;t_0)\right)\right]^{r/2}ds\right\}^{1/r} \le
  \frac{\sqrt{d}}{t_0} \left\{\int_{\bbR^d}{\frak E}(r t_0{\frak
     u}(x))\mu(dx)\right\}^{1/r}.
\end{align}
This finishes the proof of \eqref{030707}.


%
%

\subsection{Proof of Corollary \ref{Uwaga}}

\label{dowodUw}


Let 
 $
 F=fe^{-U/p}.
 $ 
Note that
\begin{align*}
&\partial^2_{x_ix_j} F=\partial_{x_ix_j}f e^{-U/p}-\frac{1}{p}\partial_{x_i}f\partial_{x_j}Ue^{-U/p}\\
&-\frac{1}{p}\partial_{x_j}f\partial_{x_i}Ue^{-U/p}-\frac1p
  f\left(\partial_{x_ix_j}U-\frac1p \partial_{x_i}U\partial_{x_j}U\right)e^{-U/p}.
\end{align*}
Let  $p'\in(p,q)$.
Using the assumption \eqref{zalozenie} for $U$ and an
appropriate $r$, we can find $C>0$ such that 
$$
\|\nabla_x^2 f\|_{L^p(\mu)}\le \|\nabla_x^2 F\|_{L^p(\bbR^d)}+C \|\nabla_x f\|_{L^{p'}(\mu)}+C \|f\|_{L^{q}(\mu)}.
$$
The classical Calderon-Zygmund $L^p$ estimates (with respect to the Lebesgue measure) give
\begin{align*}
\|\nabla_x^2 F\|_{L^p(\bbR^d)}\le C\|\Delta_x F\|_{L^p(\bbR^d)}.
\end{align*}
Note that
\begin{align*}
&\Delta_x F(x)={\cal L}f(x) e^{-U(x)/p}+a(U,\nabla_x
  U,H,\nabla_xH)\cdot \nabla_x f(x) e^{-U(x)/p}
\\
&
+b(\nabla_x U,\Delta_xU, H,\nabla_xH) f(x) e^{-U(x)/p},
\end{align*}
where $a(\cdot)$, $b(\cdot)$ are some polynomials in the indicated variables.
Therefore
$$
\|\Delta_x F\|_{L^p(\bbR^d)}\le \|{\cal L} f\|_{L^{p}(\mu)}+C \|\nabla_x f\|_{L^{p'}(\mu)}+C \|f\|_{L^{q}(\mu)}.
$$
Using the already proved estimated for $\|\nabla_x f\|_{L^{p'}(\mu)}$,
in terms of $\|{\cal L} f\|_{L^{q}(\mu)}+\|f\|_{L^{q}(\mu)}$ for $p'<q$,
we conclude the corollary.
\qed

\appendix

\section{Proofs of the auxiliary facts about the diffusion}

\subsection{Proofs of Propositions \ref{prop010804-19} and \ref{prop020804-19}}
\label{sec2.1}


 For
 $R>0$, let  $\chi_R\in C^\infty_0(\bbR^d)$ be such that $\chi_R(x)\equiv 1$ for
 $|x|\le R+1$, $\chi_R(x)\equiv 0$ for  $|x|\ge R+2$ and $\|\nabla_x\chi_R\|_\infty\le2$.
Define $H^{(R)}(x)=[H^{(R)}_{i,j}(x)]$, where
 $H^{(R)}_{i,j}(x):=H_{i,j}(x)\chi_{R}(x)$, $i,j=1,\ldots,d$ and 
 $U^{(R)}(\cdot)$ a $C^2$-smooth potential, such that 
$U^{(R)}(x)=U(x)$, $|x|\le R+1$ and $U^{(R)}(x)=|x|$, $|x|\ge R+2$.

Let $\mu_R(dx)=Z_R^{-1}e^{-U_R(x)}dx$ be a probability measure, with $Z_R$ -- the appropriate
normalizing constant. We shall assume that the truncation is made in
such a way that $\lim_{R\to+\infty}Z_R=Z$. As a result
$\lim_{R\to+\infty}\mu_R=\mu$ in total variation.

Let $L_R$ and ${\cal A}_R$ be the differential operators defined by
formulas analogous to \eqref{L} and \eqref{CA}, where the respective
coefficients have been replaced by the truncated ones introduced above. Then ${\cal
   L}_R=L_R+{\cal A}_R$ is the generator  of a family of
 diffusions 
 $\left(X_R(t,x)\right)_{t\ge0}$ indexed by the starting point. The
 stationary $X_R(t)=X_R(t,\eta)$, where $\eta(x)=x$ is
 distributed according to $\mu_R$. 
Let ${\frak e}_{R,x}:=\tau_R(X_R(\cdot,x) )$, see \eqref{tauR}, be the exit time of diffusion
$\left(X_R(t,x) \right)_{t\ge0}$ from the ball of radius $R>0$, centered at $0$. Note that
$$
X_R(t,x)=X(t,x),\quad 0\le t\le {\frak e}_{R,x}\quad\mbox{and}\quad
{\frak e}_{R,x}= \tilde {\frak e}_{R,x}:=\tau_R(X(\cdot,x) ).
$$
To prove \eqref{fe} it suffices to show that for any $T>0$
\begin{equation}
\label{031104-19}
\lim_{R\to+\infty}\int_{\bbR^d}\mu(dx)\bbP[\tilde {\frak e}_{R,x}\le T]=0.
\end{equation}
Since $\mu_R$ converges to $\mu$ in the total variation it suffices to
show that
\begin{equation}
\label{031104-19a}
\lim_{R\to+\infty}\int_{\bbR^d}\mu_R(dx)\bbP[{\frak e}_{R,x}\le T]=0.
\end{equation}
From It\^{o} formula applied to $f(X_R(t,x))$, where
$$
f(x):=(|x|^2+1)^{\rho},\quad 
$$
 we conclude
\begin{align*}
& f(X_{R}(T\wedge{\frak e}_{R,x},x))=f(x)+\frac12\int_0^{T\wedge{\frak
  e}_{R,x}}\nabla f(X_R(s,x) )\cdot\left[-\nabla_xU(X_R(s,x))+b(X_R(s,x))\right]ds\\
&+\frac12\int_0^{T\wedge{\frak
  e}_{R,x}}\Delta_x f(X_R(s,x) )ds+\int_{0}^{T\wedge{\frak e}_{R,x}}\nabla f(X_R(s,x) )\cdot  dw(s).
\end{align*}
Here, we assume that $\rho$ is as in condition A3) and it belongs to $(0,1)$.
Applying the expectation to both sides of \eqref{031104-19a} and
averaging over the initial data with respect to the  stationary
measure $\mu_R$ we conclude
\begin{align}
&\int_{\bbR^d}\mu_R(dx)\bbE f(X_{R}(T\wedge{\frak e}_{R,x},x))\le
\int_{\bbR^d}f(x)\mu_R(dx)\\
&
+\frac{T}{2}\int_{\bbR}\left|\nabla_x
  f(x)\cdot\left[-\nabla_xU(x)+b(x)\right]\right| \mu_R(dx)+2\rho Td.\nonumber
\end{align}
Here we have used the inequality
$$
|\Delta_xf(x)|=2\rho (|x|^2+1)^{\rho-2}|(d+2(\rho-1))|x|^2+d|\le
4\rho d,\quad x\in\bbR^d.
$$
Hence, there exists some contant $C>0$, independent of $R>0$, such that
$$
\int_{\bbR^d}\mu_R(dx)\bbE\left|X_{R}(T\wedge{\frak e}_{R,x},x)\right|^{2\rho}\le C+CT.
$$
As a result
$$
\int_{\bbR^d}\mu_R(dx)\bbP\left[ {\frak e}_{R,x}\le T\right]\le \frac{C+CT}{R^{2\rho}}.
$$
Passing with $R\to\infty$ we conclude \eqref{031104-19a}, which ends the
proof of Proposition \ref{prop010804-19}.\qed

\bigskip

To show 
Proposition \ref{prop020804-19}, observe first that an analogue of
\eqref{stat} holds for $\left(X_R(t)\right)_{t\ge0}$ and functions
$f_1,\ldots,f_n$ that are bounded and continuous.
Letting $R\to+\infty$ we conclude the equality for stationary process $\left(X(t)\right)_{t\ge0}$
and bounded continuous functions. Using an approximation argument we
can extend \eqref{stat} to arbitary bounded and measurable functions
$f_1,\ldots,f_n$ that ends the proof of Proposition
\ref{prop020804-19}. \qed

\subsection{Proof of Proposition \ref{prop030804-19}}

\label{sec2.3}

From Proposition \ref{prop010804-19} it follows that
\begin{equation}
\label{031104-19b}
\int_{\bbR^d}\mu(dx)\bbP[\tilde {\frak e}_{x}<+\infty]=0.
\end{equation}
The above equality implies that there exists a Borel measurable set
of null  Lebesgue measure   ${\cal Z}\subset \bbR^d$,
such that for any $x\not\in {\cal Z}$ we have $\tilde {\frak
  e}_{x}=+\infty$.

For an arbitrary $x\in\bbR^d$ 
there exists $R>0$ such that the sphere $S_R(x):=[y:\,|y-x|=R]$
intersects with $ {\cal Z}$ on a
set ${\cal N}$ of a null surface Lebesgue measure $\om_{d,R}$. Consider the
harmonic measure 
$\om_R^x(A)=\bbP[X_R({\frak e}_{R,x},x)\in A]$, where $A\in{\cal B}(S_R(x))$ --
the Borel $\si$-algebra of subsets of the sphere. It is equivalent
with respect to $\om_{d,R}$, see Theorem 4.4, p. 311 of
\cite{pinsky}. Hence  $\om_R^x({\cal N})=0$. As a result 
$$
X_R({\frak e}_{R,x},x)=X(\tilde {\frak e}_{R,x},x)\not \in {\cal
  N},\quad \mbox{a.s.}
$$
and for such an event we can let
$$
X(t,x):=X(t-\tilde {\frak e}_{R,x} ,X(\tilde {\frak e}_{R,x}))
\quad \mbox{for all $t\ge \tilde {\frak e}_{R,x}$},
$$
which ends the proof of  Proposition \ref{prop030804-19}.\qed

\subsection{Proof of Proposition \ref{prop040804-19}}

\label{sec2.4}

We use the notation of Section \ref{sec2.1}.
Let $\nu$ be the law of the random vector $(X(0,\eta),X(1,\eta))$ in
$(\bbR^{2d},{\cal B}(\bbR^{2d}))$. We claim that it is equivalent with
the Lebesgue measure $m_{2d}$, i.e. the families of  null sets in both
$\nu$ and $m_{2d}$ measures are equal. First suppose that $A\in {\cal
  B}(\bbR^{2d})$ is a $m_{2d}$-null measure set. Given $R>0$ let
$\nu_R$ be the
law of $(X_R(0,\eta),X_R(1,\eta))$, with $\eta$ distributed according
to $\mu_R$. It is quite straightforward (e.g. using the
Girsanov theorem) to conclude that $\nu_R$ is equivalent with
$m_{2d}$.
Let $A_{R}:=A\cap(S_R(0)\times S_R(0))$ and  let ${\frak e}_{R,x}:=\tau_R(X_R(t,x) )=\tau_R(X(t,x) )$.  
We have $\nu_R[A_R]=0$ and, by the argument made in
Section \ref{sec2.1}, we know  that
 \begin{equation}
\label{011304-19}
\lim_{R\to+\infty}\int_{\bbR^d}\,\mu(dx)\bbP[{\frak
  e}_{R,x}< 1]=0.
\end{equation}
Hence
\begin{align*}
&\nu [A_{R}]=\int_{\bbR^d}\,\mu(dx)\bbP[(x,X(1,x))\in A_{R}]\\
&
\le
  \int_{\bbR^d}\,\mu(dx)\bbP[(x,X(1,x))\in A_{R},\,{\frak
  e}_{R,x}\ge 1]+ \int_{\bbR^d}\,\mu(dx)\bbP[{\frak
  e}_{R,x}< 1]\\
&
\le \frac{Z_R}{Z}\nu_R[ A_{R}]+ \int_{\bbR^d}\,\mu(dx)\bbP[{\frak
  e}_{R,x}< 1]=\int_{\bbR^d}\,\mu(dx)\bbP[{\frak
  e}_{R,x}< 1].
\end{align*}
Letting $R\to+\infty$, we conclude that $\nu [A]=\lim_{R\to+\infty}\nu
[A_{R}]=0$.

Now suppose that  $m_{2d}(A)>0$ and   $\nu [A]=0$. Then, for a
sufficiently large $R>0$ we would have $m_{2d}(A_R)>0$ and 
\begin{align*}
&\nu [A]\ge\nu [A_{R}]\ge 
  \int_{\bbR^d}\,\mu(dx)\bbP[(x,X(1,x))\in A_{R},\,{\frak
  e}_{R,x}\ge 1]\\
&
=\frac{Z_R}{Z}\int_{\bbR^d}\,\mu_R(dx)\bbP[(x,X_R(1,x))\in A_R ,\,{\frak
  e}_{R,x}\ge 1]>0.
\end{align*}
The last strict inequality can be seen  using, for example, the Girsanov
theorem.

Using the Radon-Nikodym theorem we conclude therefore the existence of
a  strictly positive, Borel measurable density $p(x,y)$,  such that 
$$
\bbE_{\mu} F(X(0),X(1))=\int_{\bbR^{2d}}F(x,y)p(x,y)\mu(dx) \mu(dy)
$$
for any Borel measurable and bounded function
$F:\bbR^{2d}\to\bbR$. Here $\bbE_{\mu}$ is the expectation with
respect to measure $\mu\otimes\bbP$. This in particular implies
\eqref{polgrupa}.

\end{document}